\documentclass[journal]{IEEEtran}
\usepackage{pgfplots}
\usepackage{dblfloatfix} 

\usepackage{mathtools, cuted}
\usepackage{mynotes}
\usepackage{cancel}
\definecolor{amethyst}{rgb}{1, 0, 1}
\definecolor{blue-violet}{rgb}{0.54, 0.17, 0.89}
\definecolor{brightturquoise}{rgb}{0.03, 0.91, 0.87}
\usepackage{pgfplots}
\usepackage{caption}
\usepackage{multirow}

\graphicspath{{./figs/}}

\title{ Efficient Robust Parameter Identification in Generalized Kalman Smoothing Models}

 \author{%
Jonathan Jonker$^*$\thanks{$^*$Department of Mathematics, University of Washington, Seattle, WA (jonkerjo@uw.edu)}
\and    Peng Zheng$^\dagger$\thanks{$^\dagger$Department of Applied Mathematics, University of Washington, Seattle, WA
(zhengp@uw.edu)} 
\and    Aleksandr Y. Aravkin$^\dagger$\thanks{$^\dagger$Department of Applied Mathematics, University of Washington, Seattle, WA
(saravkin@uw.edu)}}

\begin{document}

\maketitle

\begin{abstract}
Dynamic inference problems in autoregressive (AR/ARMA/ARIMA), exponential smoothing, and navigation  
are often formulated and solved  using state-space models (SSM), which  
allow a range of statistical distributions to inform innovations and errors. 
In many applications the main goal is to identify not only the hidden state, but also 
additional unknown model parameters (e.g. AR coefficients or unknown dynamics).

We show how to efficiently optimize over model parameters in SSM that use smooth process and measurement losses.  
Our approach is to project out state variables, obtaining a value function that only depends on the parameters of interest, 
and derive analytical formulas for first and second derivatives that can be used by many types of optimization methods. 
%We draw a distinction between nonsingular and singular formulations, since 
% most of the models of interest have singular covariance structure. 

The approach can be used with smooth robust penalties such as Hybrid and the Student's t, in addition to classic least squares. 
We use the approach to estimate robust AR models and long-run unemployment rates with sudden changes.  
 \end{abstract}

\section{Introduction.}

The linear state space model is widely used in tracking and navigation~\cite{ybarshalom-2001a}, control~\cite{anderson2007optimal},  signal processing \cite{AndersonMoore}, and  other time series~\cite{hyndman2002state,tsay2005analysis}. 
The model assumes linear relationships between latent states with noisy observations: 
\begin{equation} \label{eq:statespace}
\begin{aligned}
x_k & =  G_kx_{k-1} + \epsilon_k^p, \quad  && k = 1, \dots, N, \\
z_k & =  H_kx_k + \epsilon_k^m, \quad  && k = 1, \dots, N,
\end{aligned}
\end{equation}
where $x_0$ is a given initial state estimate, $x_1, \dots, x_N$ are unknown latent states with known linear process models $G_k$, and $z_1, \dots, z_N$ are observations obtained using known linear models $H_k$. 
The errors $\epsilon_k^p$ and $\epsilon_k^m$  are assumed to be mutually independent random variables with covariances $Q_k$ and $R_k$.
These covariances may be singular to capture standard autoregressive structures.

In many applications the models $G_k, H_k, Q_k, R_k$ are specified up to model parameters $\theta$. 
We restrict out attention to formulations where variances $Q_k, R_k$ are known, 
while $G_k(\theta)$ and $H_k(\theta)$ are  $\mathcal{C}^2$ mappings of $\theta$. 
This captures smoothing parameters in Holt-Winters c.f.~\cite{hyndman2002state}, 
autoregressive and moving average parameters in ARMA c.f.~\cite{tsay2005analysis}, 
and unknown dynamic parameters in navigation models. In most of these models, 
$G$ and $H$ are affine functions of unknown parameters $\theta$.

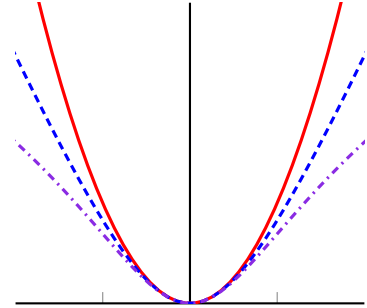
\begin{figure}[h!]
\center
\begin{tikzpicture}
  \begin{axis}[
    thick,
    height=4cm,
    xmin=-2,xmax=2,ymin=0,ymax=1.5,
    no markers,
    samples=50,
    axis lines*=left, 
    axis lines*=middle, 
    scale only axis,
    xtick={-1,1},
    xticklabels={},
    ytick={0},
    ] 
%\addplot[red,domain=-2:-1,densely dashed]{-x-.5};
\addplot[red, domain=-2.5:+2.5, very thick]{.5*x^2};
\addplot[blue, domain=-2.5:+2.5, densely dashed, very thick]{sqrt(1+(x/1)^2) - 1};
\addplot[blue-violet, domain=-2.5:+2.5, dash dot, very thick]{0.5*ln(1+x^2)};
%\addplot[red,domain=+1:+2,densely dashed]{x-.5};
%\addplot[blue,mark=*,only marks] coordinates {(-1,.5) (1,.5)};
  \end{axis}
  \end{tikzpicture}
\caption{\label{fig:losses} Common smooth loss functions: least-squares (red solid), Hybrid (blue dashed), and student's T (violet dash dot).}
\end{figure}

Standard models assume the errors $\epsilon_k^p$ and $\epsilon_k^m$ are Gaussian, which gives rise to 
the least squares penalty in the inference problem, see the red solid curve in Figure~\ref{fig:losses}.
Changing the observation model to the Hybrid (blue dash) or Student's t loss (violet dash dot) 
 robustifies model estimates in the face of measurement outliers. Analogous changes to the innovations model
allows the framework to track sudden changes.

{\bf Motivating Application: Structural Unemployment Rate}. 
We are interested in fitting parameters within structural unemployment rate models, see e.g.~\cite{mocan1999structural}. 
The state vector 
\begin{equation}
\label{eq:state}
x_k = \begin{bmatrix}
u_{k-1} & u_{k-1}^c & u_k & u_k^c
\end{bmatrix}^T\\
\end{equation}
 tracks total ($u$) and `cyclic' ($u^c$) unemployment using 
the auto-regressive model 
\begin{equation}
\label{eq:dynamics}
% x_k = \begin{bmatrix}
%u_{t-1}^*\\
%u_{t-1}^c\\
%r_k^p^*\\
%r_k^p^c
%\end{bmatrix}, \quad
 G_k = \begin{bmatrix}
0 & 0 & 1 & 0\\
0 & 0 & 0 & 1\\
-l_1 & 0 & 1+l_1 & 0\\
0 & 1/2-l_2 & 0 & l_2
\end{bmatrix}, 
\quad \epsilon_k^p = \begin{bmatrix}
0\\
0\\
\epsilon_k^1\\
\epsilon_k^2
\end{bmatrix},
\end{equation}
%Define our observation equation as,
%
%
%
\begin{equation}
\label{eq:measurements}
% z_k = \begin{bmatrix}
%r_k^p\\
%\Delta w_k^c
%\end{bmatrix}, \quad
 H_k = \begin{bmatrix}
0 & 0 & 1 & 1\\
0 & \gamma/2 & 0 & \gamma/2
\end{bmatrix}, \quad
 \epsilon_k^m = \begin{bmatrix}
\epsilon_k^3\\
\epsilon_k^4
\end{bmatrix}.
\end{equation}
Here, $\ell_1$ and $\ell_2$ are auto-regressive parameters 
while $\gamma$ is an unknown measurement parameter.  
Unemployment rates can experience fast changes, so 
we need a heavy tailed model for innovations.   
To solve the full problem, we must 

\begin{enumerate}
\item Estimate  states $\{x_k\}$ as well as parameters $\ell_1, \ell_2$, $\gamma$.
\item Account for the singular process covariance $Q$.
\item Use non-Gaussian losses (e.g. Hybrid and Student's t) for $\epsilon_k^p$ to track fast rate changes.
\end{enumerate}

The paper proceeds as follows. In Section~\ref{sec:theory},
%~\ref{sec:formulations}, 
we review optimization formulations for state and model parameter 
inference using singular and nonsingular covariance models, 
and introduce the {\it value function} which depends only on the model parameters, 
as e.g. $\ell_1, \ell_2$, $\gamma$ above.  
In Sections~\ref{sec:nonsingular_ssm} and~\ref{sec:singular_ssm}, we look in detail at nonsingular and singular 
Kalman smoothing models, and obtain existence results 
and formulas for first and second derivatives of the value function.  
  Finally, in Section~\ref{sec:numerics}
we present use cases that show how to efficiently obtain structural parameters
when using general losses and singular covariance structure. 

%%%%%
\section{Notations and Preliminaries}
\label{sec:notations}
We first introduce notation and  key definitions. %that are used in the paper.
\paragraph{Superscript and subscript}
We use superscripts to distinguish process ($p$) and measurement ($m$) model variables and subscripts
to represent partial derivatives ($\theta, y, r,\ldots$) and the index in the Kalman model ($k$).
When taking derivatives, subscripts indicate the position rather than actual variable.
% Subscript letters identify the variable with respect to which we take partial derivative.
% of the partial derivative rather than that of the variables.
For example, if we denote $v(\theta) = f(\theta, Y(\theta))$, then 
\[
\nabla v(\theta) \ne f_\theta(\theta, Y(\theta)) = \partial_\theta f(\theta, y)|_{y = Y(\theta)}.% \ne \nabla v(\theta).
\]
\paragraph{Loss functions}
We use the following loss functions:
\begin{itemize}
\setlength\itemsep{1em}
\item Least squares:
$
\ell(r) = \frac{1}{2}\|r\|^2.
$
\item Hybrid:
$
\ell(r; \nu) = \sum_i \sqrt{r_i^2 + \nu^2} - \nu.
$
\item Student's T:
$
\ell(r; \nu) = \sum_i \ln(1 + r_i^2/\nu).
$
\end{itemize}

%%%%%
\section{Differentiating Implicit Functions}
\label{sec:theory}
In this section, we introduce  a general theoretical result for calculating the derivatives for implicit functions
in an optimization context.
We then specialize this general theorem to nonsingular and singular state space models (SSM) in the next section.
Consider a $\cC^2$-smooth function, $f: \bR^n \times \bR^m \rightarrow \bR$,
where  in SSM models we denote parameters by $\theta$ in $f(\theta, y)$ and the states 
together with any auxiliary variables (such as dual variables) by $y$.
The appropriate stationary condition is given by 
\begin{equation}
\label{eq:stationary}
\cH(\theta, y) := f_y(\theta, y) = 0.
\end{equation}
For any given $\theta$, the optimal estimate $y(\theta)$ 
is obtained by solving the equation $\cH(\theta, y) = 0$;
in particular $y(\theta)$ depends on $\theta$ implicitly. 
%which give us implicit dependence of $y$ on $\theta$.
%
We introduce a variant of the implicit function theorem presented by~\cite{bell2008algorithmic} %to study the property of $\cH$.
to characterize the structure of this implicit dependence. 
\begin{theorem}[Implicit Functions and Derivatives]
\label{thm:implicit}
Suppose that $\mathscr{U} \subset\mathbb{R}^n$ and $\mathscr{V} \subset \mathbb{R}^m$ are open, 
$\cH: \mathscr{U} \times \mathscr{V} \rightarrow \mathbb{R}^m$ is continuously differentiable.
If there exists $\overline \theta \in \mathscr{U}$ and $\overline y \in \mathscr{V}$, such that $\cH(\overline \theta, \overline y) = 0$
and $\cH_y(\overline \theta, \overline y)$ is invertible.
Then there exists (if necessary we choose $\mathscr{U}$ and $\mathscr{V}$ to be small neighborhood of $\overline \theta$ and $\overline y$ to guarantee the existence) a $\cC^1$ mapping $Y:\mathscr{U} \rightarrow \mathscr{V}$ satisfying $Y(\overline \theta) = \overline y$, 
and $\cH(\theta, Y(\theta)) = 0$ for all $\theta$ in $\mathscr{V}$. 

Moreover, we have the formula 
\[
Y_\theta(\theta) = - \cH_y(\theta, Y(\theta))^{-1}\cH_\theta(\theta, Y(\theta)).
\]
\end{theorem}
When the function $Y(\theta)$ as above exists, we can 
define the value function %under the assumption of Theorem~\ref{thm:implicit},
\begin{equation}
\label{eq:value}
v(\theta) = f(\theta, Y(\theta)).
\end{equation}
Our goal is to compute first and second derivatives of $v$,
which are summarized in the following corollary.
\begin{corollary}[Derivatives of the Value Function~\eqref{eq:value}]
\label{cor:val_der}
Under the assumptions of Theorem~\ref{thm:implicit}, 
and using $\overline y$ to represent the $y$ obtained by evaluating 
$Y(\overline \theta)$, we have,
\begin{equation}
\label{eq:derivatives}
\begin{aligned}
v_\theta(\overline \theta) =& f_\theta(\overline \theta, \overline y)\\
v_{\theta\theta}(\overline \theta) =& f_{\theta\theta}(\overline \theta, \overline y)- \cH_{\theta}(\overline \theta, \overline y)^\top \cH_y(\overline \theta, \overline y)^{-1} \cH_{\theta}(\overline \theta, \overline y).
\end{aligned}
\end{equation}
\end{corollary}
These derivations are along the lines of those presented by Bell and Burke~\cite{bell2008algorithmic} and are mainly given here for a self-contained exposition.
We now compute analytic expressions of derivatives with respect to model parameters 
for both nonsingular and singular Kalman smoothing systems.

%%%%%
\section{Nonsingular SSM}
\label{sec:nonsingular_ssm}
Consider the case where the covariance matrices $Q_k$ and $R_k$ in SSM are nonsingular.
Pre-whitening $\epsilon_k^p$ and $\epsilon_k^m$, the objective function of interest is given by 
\begin{equation}
\label{eq:nonsingular}
\begin{aligned}
f(\theta,y) = \sum_{k=1}^N  &\lt\{\ell_k^p\lt(Q_k^{-1/2}(x_k - G_k(\theta) x_{k-1}) \rt)\rt. \\
+ &\lt.\ell_k^m\lt(R_k^{-1/2} (z_k - H_k(\theta) x_k)\rt)\rt\},
\end{aligned}
\end{equation}
where $y = x = [x_1;\ldots;x_N]$, $\ell_k^p$ and $\ell_k^m$ are the loss function corresponding to the distributions of
$\epsilon_k^p$ and $\epsilon_k^m$.
Here we assume $\ell_k^p$, $\ell_k^m$ are smooth; three key examples 
are least squares, Hybrid, and Student's t losses introduced in Section~\ref{sec:notations}. 
Objective \eqref{eq:nonsingular} can be written compactly as,
\begin{equation}
\label{eq:compact_n}
\begin{aligned}
f(\theta, y) =\,&\ell^p\lt(Q^{-1/2}(G(\theta)x - \zeta)\rt) + \\
&\ell^m\lt(R^{-1/2}(H(\theta)x - z)\rt),
\end{aligned}
\end{equation}
where,
\begin{align*}
&G(\theta)   = \begin{bmatrix}
    I  & 0      &          &
    \\
    -G_2(\theta)   & I  & \ddots   &
    \\
        & \ddots &  \ddots  & 0
    \\
        &        &   -G_N(\theta)  & {I}
\end{bmatrix},
&&Q = \begin{bmatrix}
Q_1 & &\\
& \ddots &\\
& & Q_N
\end{bmatrix},\\
&H(\theta) = \begin{bmatrix}
    H_1(\theta)  & 0      &          &
    \\
    0   & H_2(\theta)  & \ddots   &
    \\
        & \ddots &  \ddots  & 0
    \\
        &        &   0  & H_N(\theta)
\end{bmatrix},
&&R = \begin{bmatrix}
R_1 & &\\
& \ddots &\\
& & R_N
\end{bmatrix},
\end{align*}
and $\zeta = [x_0; 0; \ldots; 0]$ and $z = [z_1; \ldots; z_N]$.
The stationary condition in this case is given by 
\begin{equation}
\label{eq:h_nonsingular}
\begin{aligned}
\cH(\theta, y) &= f_y(\theta, y)\\
&= G(\theta)^\top Q^{-\top/2}\ell_r^p(r^p) + H(\theta)^\top R^{-\top/2}\ell_r^m(r^m)
\end{aligned}
\end{equation}
where $r^p = Q^{-1/2}(G(\theta)y - \zeta)$ and $r^m = R^{-1/2}(H(\theta)y - z)$.
%%%%%

In the least squares case,  \eqref{eq:h_nonsingular} is a linear equation  solved 
by inverting the block-tridiagonal system 
\[
G(\theta)^\top Q^{-1} G(\theta) + H(\theta)^\top R^{-1} H(\theta). 
\]
for more general smooth penalties $\ell^p, \ell^m$ a Newton method 
is needed to compute $y(\theta)$. 

By Theorem~\ref{thm:implicit}, existence and differentiability 
of  $Y(\theta)$ is guaranteed by the existence of the pair $(\overline \theta, \overline y)$ such that,
$\cH(\overline \theta, \overline y) = 0$ and the partial Hessian below is invertible: 
\[\begin{aligned}
\cH_y(\overline \theta, \overline y) = 
&G(\overline \theta)^\top  Q^{-\top/2}\ell_{rr}^p(r^p) Q^{-1/2}G(\overline \theta) + \\
&H(\overline \theta)^\top  R^{-\top/2}\ell_{rr}^m(r^m) R^{-1/2}H(\overline \theta).
\end{aligned}\]
When $\ell^p, \ell^m$ are least squares or Hybrid, $\cH_y( \theta,  y)$ is invertible for any $(\theta, y)$,
and for every $\overline \theta$ there exist a $\overline y$ such that $\cH(\overline \theta, \overline y) = 0$,
since both penalties are strictly convex. In the case of the Student's t, the Hessian 
may fail to be positive definite at some pairs $(\theta, y)$~\cite{aravkin2014robust} and there is no absolute guarantee
that the methodology will hold, as expected for a potentially nonconvex formulation.
Practical behavior is another matter, and in numerical experience we see a positive definite Hessian at the minimizer, 
and so the derivative formulas hold. A safeguard can be added to any practical implementation, 
that can trigger a failsafe and use a less efficient method.
%so the derivative formulas hold.   

We now compute remaining terms in Corollary~\ref{cor:val_der}, assuming for simplicity 
that $G(\theta)$ and $H(\theta)$ are affine functions of $\theta$. 
\[
\ol{r^p}  := Q^{-\top/2}\ell_r^p(r^p), \quad \ol{r^m} := R^{-\top/2}\ell_r^m(r^m)
\]
\[
\begin{aligned}
f_\theta(\theta, \overline y) &= (G \ol{x})_\theta^\top  Q^{-\top/2}\ell_r^p(r^p) + (H\ol{x})_\theta^\top  R^{-\top/2}\ell_r^m(r^m)\\
f_{\theta\theta}(\theta, \ol{y}) & = (G \ol{x})_\theta^\top  Q^{-\top/2}\ell_{rr}^p(r^p)Q^{-1/2}(G \ol{x})_\theta \\ 
& + (H \ol{x})_\theta^\top  R^{-\top/2}\ell_{rr}^m(r^m)R^{-1/2}(H \ol{x})_\theta\\
\cH_\theta (\theta, \ol{y})& =  (G(\theta)^\top \ol{r^p})_\theta + G(\theta)^\top Q^{-\top/2}\ell_{rr}^p(r^p) Q^{-1/2}(G\overline x)_\theta \\
& + (H(\theta)^\top \ol{r^m})_\theta + H(\theta)^\top R^{-\top/2}\ell_{rr}^m(r^m) R^{-1/2}(H\overline x)_\theta
\end{aligned}
\]
We now have, fully and explicitly, first and second derivatives of the value function $v(\theta)$ in~\eqref{eq:derivatives}
for the nonsingular case. 
Though these results are straightforward, they do not appear in any smoothing literature we are aware of 
in  this compact form, even for least squares losses.

\section{Singular SSM}
\label{sec:singular_ssm}

When covariances $Q_k$ and $R_k$ are singular, 
we rewrite~\eqref{eq:nonsingular} to include null-space 
constraints. A singular covariance matrix 
precludes any errors and innovations that are not in its range. 
We follow~\cite{jonker2019fast}  in formulating this problem: 
\begin{equation}
\label{eq:sing}
\begin{aligned}
\min_{\theta, x, r^p, r^m}~~&\ell^p (r^p)  + \ell^m(r^m)\\
\mbox{s.t.}~~& Q^{1/2}r^p = G(\theta)x - \zeta,\\
&R^{1/2} r^m = H(\theta)x - z,
\end{aligned}
\end{equation}
and introduce the {\it Lagrangian}
\begin{equation}
\label{eq:Lag}
\begin{aligned}
f(\theta, y) =& \cL(\theta, x, r^p, r^m, \lambda^p, \lambda^m) \\
=& \ell^p (r^p)  + \ell^m(r^m) \\
&- \ip{\lambda^p, Q^{1/2}r^p - G(\theta)x + \zeta } \\
&-\ip{\lambda^m, R^{1/2}r^m - H(\theta)x + z }.
\end{aligned}
\end{equation}

When $Q_k$ and $R_k$ are invertible, we can solve for $r^p, r^m$ in~\eqref{eq:sing} 
and reduce the problem to~\eqref{eq:nonsingular},
%Objective~\eqref{eq:sing} can also be written using range
%constraints: 
%\begin{equation}
%\label{eq:sing2}
%\begin{aligned}
%&\sum_{k=1}^N  \ell^p\left(Q_k^{\dagger/2}(x_k - G_k x_{k-1}) \right) + \ell^m\left(R_k^{\dagger/2} (z_k - H_k x_k)\right), \\
%& \mbox{s.t.} \quad z_k - H_k x_k \in \mathcal{R}(R_k), \quad x_k - G_k x_{k-1} \in \mathcal{R}(Q_k).
%\end{aligned}
%\end{equation}
%The $r_k^p$ and $r_k^m$ variables in~\eqref{eq:sing} give us particular 
%vectors that realize the range constraints in~\eqref{eq:sing2}. 
%While~\eqref{eq:sing2} is useful for understanding the nature 
%of singular formulations, 
so nonsingular systems are a special case of~\eqref{eq:sing}. 
Formulation~\eqref{eq:sing} can be solved
for a variety of loss functions $\ell^p$ and $\ell^m$~(see \cite{jonker2019fast}). 

Here we define the {\it value function} as a mini-max problem using the Lagrangian:  
\begin{equation}
\label{eq:value_lag}
\begin{aligned}
v(\theta) &:= \max_{\lambda^p, \lambda^m} \min_{x, r^p, r^m} \cL(\theta, x, r^p, r^m, \lambda^p, \lambda^m). 
\end{aligned}
\end{equation}
The system of equations we are interested in is now 
\[
0 = \cH(\theta, y) := f_y(\theta, y), \quad y := \{x, r^p, r^m, \lambda^p, \lambda^m\};
\]
that is, the system of equations that defines a saddle point of the Lagrangian. 
Explicitly,   $f_y(\theta, y) = 0$ is given by 
\[
f_y(\theta, y) = 
\begin{bmatrix}
G(\theta)^\top \lambda^p + H(\theta)^\top \lambda^m\\
\ell_r^p(r^p) - Q^{\top/2}\lambda^p\\
\ell_r^m(r^m) - R^{\top/2}\lambda^m\\
G(\theta)x - \zeta - Q^{1/2} r^p\\
H(\theta)x - z - R^{1/2} r^m
\end{bmatrix} = 0
\]
$Y(\theta)$ is differentiable when $\cH_y(\theta, y) = f_{yy}(\theta, y)$ is invertible:
\begin{equation}
\label{eq:HessLag}
f_{yy}(\theta, y)= \begin{bmatrix}
0 & 0 & 0 & G(\theta)^\top & H(\theta)^\top \\
0 & \ell_{rr}^p(r^p) & 0 & -Q^{\top/2} & 0 \\
0 & 0 & \ell_{rr}^m(r^m)  & 0 & -R^{\top/2} \\
G(\theta) & -Q^{1/2} & 0 & 0 & 0 \\
H(\theta) & 0 & -R^{1/2} & 0 & 0 
\end{bmatrix}
\end{equation}

We state the following theorem. 
\begin{theorem}
\label{thm:invert}
Invertibility of $f_{yy}(\theta, y)$ is equivalent to invertibility of the so called Hessian of the Lagrangian
\begin{equation}
\label{eq:SchurLag}
\begin{aligned}
&H(\theta)G(\theta)^{-1}Q^{1/2} (\ell_{rr}^p(r^p))^{-1}Q^{\top/2}G(\theta)^{-\top}H(\theta)^\top   \\
&+ R^{1/2} (\ell_{rr}^m(r^m))^{-1}R^{\top/2}.
\end{aligned}
\end{equation}
When $\ell^p, \ell^m$ are the least squares or Hybrid penalties,~\eqref{eq:SchurLag} is invertible if and only if 
\begin{equation}
\label{eq:NullCond}
\cN(R) \cap \cN(QG^{-T}H^\top ) = \{0\}.
\end{equation}

\end{theorem}

%We immediately obtain a simple corollary that is applicable to the motivating application~\eqref{eq:model}
%\begin{corollary}
%$\cX(\theta)$ is differentiable when $R$ is nonsingular.
%\end{corollary}
When Theorem~\ref{thm:invert} holds, we use Corollary~\eqref{cor:val_der} to get derivatives
of $v(\theta)$ in~\eqref{eq:value_lag}:
\begin{equation}
\label{eq:singular_ssm_der}
\begin{aligned}
v(\overline \theta) =& f(\overline \theta, \overline y)\\
v_\theta(\overline \theta) =& f_\theta(\overline \theta, \overline y)\\
v_{\theta\theta}(\overline \theta) =& f_{\theta\theta}(\overline \theta, \overline y)- \cH_{\theta}(\overline \theta, \overline y)^\top \cH_y(\overline \theta, \overline y)^{-1} \cH_{\theta}(\overline \theta, \overline y).
\end{aligned}
\end{equation}
It remains only to compute $f_\theta$, $f_{\theta\theta}$, $\cH_{\theta}$, and $\cH_{y}$.
\[
f_\theta(\theta, \ol{y}) := (\ip{\ol{\lambda^p}, G(\theta)\ol{x}})_\theta + (\ip{\ol{\lambda^m}, H(\theta)\ol{x}})_\theta
\]
When $G, H$ are affine functions of $\theta$, we have $f_{\theta\theta} = 0$. 
Finally, 
\[
\cH_{\theta}(\theta, \ol{y}) = f_{y\theta}(\theta, \ol{y}) = 
\begin{bmatrix}
(G(\theta)^\top  \ol{\lambda^p})_\theta +  (H(\theta)^\top \ol{\lambda^m})_{\theta} \\
0\\
0\\
(G(\theta)\ol{x})_\theta \\
(H(\theta)\ol{x})_\theta
\end{bmatrix}.
\]

\subsection{Special case: Invertible $R$.}

The structural unemployment model in the introduction has a singular $Q$ but an invertible $R$. 
In such cases, the derivative formulas can be written 
using only primal quantities, which significantly decreases the notational burdern. 
In particular, using the optimality conditions we have 
\[
\begin{aligned}
\lambda^m &= R^{-\top/2}\ell_r^m(r^m)\left(R^{-1/2}(H(\theta)y-z)\right)\\
\lambda^p & = - G(\theta)^{-\top}H(\theta)^\top \lambda^m
\end{aligned}
\]
and so we get the explicit primal-only formula for $v_\theta(\theta)$ 
by plugging these expressions into 
\[
v_\theta(\theta) =  (\ip{\ol{\lambda^p}, G(\theta)\ol{y}})_\theta + (\ip{\ol{\lambda^m}, H(\theta)\ol{y}})_\theta.
\] 

\subsection{Special case: Least-squares.}

If $\ell^p(\cdot)$ and $\ell^m(\cdot)$ are both given by $\frac{1}{2}\|\cdot\|^2$, 
the optimality conditions simplify substantially, and we have 
\[
\begin{aligned}
r^p & =  Q^{T/2}\lambda^p  , \quad r^m = R^{T/2}\lambda^m\\
Q^{1/2}r^p &= Q\lambda_p  =   G(\theta) x - \zeta, \\
 R^{1/2}r^m & = R\lambda_m  = H(\theta) x - z \\
0 & = G(\theta)^\top \lambda^p + H(\theta)^\top \lambda^m.
\end{aligned}
\]
Plugging these conditions back into the Lagrangian, we get the dual objective 
\[
\begin{aligned}
f(\lambda^p, \lambda^m) & = -\frac{1}{2}(\lambda^p)^\top Q\lambda^p  - \frac{1}{2}(\lambda^m)^\top R\lambda^m  - (\lambda^p)^\top  \zeta - (\lambda^m)^\top z \\
\mbox{s.t.}\quad  &  G(\theta)^\top \lambda^p + H(\theta)^\top \lambda^m = 0.
\end{aligned}
\]
Using invertibility of $G$, we eliminate $\lambda_p$: 
\[
\begin{aligned}
f(\lambda_m) &= -\frac{1}{2}(\lambda^m)^\top \left( HG^{-1}QG^{-T}H^\top  + R\right)\lambda^m\\ 
&- (\lambda^m)^\top (z - HG^{-1}\zeta).
\end{aligned}
\]
In the least squares case, the dual solution $\lambda^m$ is unique exactly when the linear system 
\begin{equation}
\label{eq:dual_inv}
HG^{-1}QG^{-T}H^\top  + R
\end{equation}
is invertible, and then we have 
\[
\lambda_m  = \left( HG^{-1}QG^{-T}H^\top  + R\right)^{-1}(z - HG^{-1}\zeta),
\]
a closed form solution. 
A simple sufficient condition for the invertibility of~\eqref{eq:dual_inv} is to have $R$ itself invertible,
as in the special case previously discussed. 

%which can be plugged into~\eqref{eq:derSingular} to get a complete closed form solution. 
%However, evaluating this formula requires evaluating a dense system. 
%We get the dual solution $\lambda_m$ for free when we use techniques of~\cite{jonker2019fast}
%to solve the singular smoothing system; and we need not rely on the least squares case. 

\section{Numerical Examples}
\label{sec:numerics}

We now apply the results of the previous sections to analyze two simple singular models with unknown 
states and parameters.  In Section~\ref{sec:ar} we present a state-space formulation 
for the AR-1 model, show how to robustify it to outliers in the data, and present explicit derivatives for the value function. 
We use these derivatives to design an efficient solver for both standard and robust AR models. 
In Section~\ref{sec:rates}, we apply the methods in this paper to fit a structural model for 
unemployment rates that can track fast changes. While the structural unemployment model 
is currently used in the EU, in this paper we only show  results on simulated synthetic data 
where we know ground truth, and leave any more detailed work to future collaborations.

\subsection{Robust AR Fitting}
\label{sec:ar}
An AR-1 model begins with equations
\begin{equation}
\begin{aligned}
x_k &= c + \varphi x_{k-1} + \epsilon_k^p\\
y_k &= H_k x_k + \epsilon_k^m
\label{eq:AR1}
\end{aligned}
\end{equation}
Where $\epsilon_k^p, \epsilon_k^m$ have covariances $Q_k, R_k$ respectivly, $c$ is an unknown constant, and $\varphi$ is a parameter to be estimated. To make this take the form of~(\ref{eq:statespace}) we create an augmented state
\begin{equation}
\hat{x}^k = \begin{bmatrix}x_k\\ c\end{bmatrix}
\label{eq:aug}
\end{equation}

and use state equations
\begin{equation}
\begin{aligned}
\label{eq:augstate}
\hat{x}^k = \begin{bmatrix}\varphi & 1\\ 0 & 1\end{bmatrix} \hat{x}_{k-1} + \hat{\epsilon}_k^p\\
y_k = \begin{bmatrix} H_k & 0\\ 0 & 0\end{bmatrix}\hat{x}_k + \hat{\epsilon}_k^m
\end{aligned}
\end{equation}
Then
\[
\hat{Q}_0 = \begin{bmatrix}Q_0 & 0\\ 0 & 1\end{bmatrix}, \quad \hat{Q}_k = \begin{bmatrix}Q_k & 0\\ 0 & 0\end{bmatrix} \quad k > 0
\]
\[
\hat{R}_k = \begin{bmatrix}R_k & 0\\ 0 & 0\end{bmatrix}
\]
This choice of $\hat{Q}_0$ will allow us to fit the constant $c$ as part of the state while $\hat{Q}_k, k>0$ will act as equality constraints holding it constant through all time points. In order to compute the derivatives of the value functions we first note the following derivative formula in this case for any vector $\eta$.
\begin{equation}
\label{eq:ARder}
\Big(G_i(\theta)\eta\Big)_\theta = \Big(\begin{bmatrix}\varphi & 1\\ 0 & 1\end{bmatrix}\begin{bmatrix}\eta^1 \\ \eta^2\end{bmatrix}\Big)_\theta = \begin{bmatrix}\eta^1 \\ 0\end{bmatrix} = \tilde D\eta
\end{equation}
Where $\tilde D = \begin{bmatrix}1 & 0\\ 0 & 0\end{bmatrix}$.
Define
\[
D = \begin{bmatrix}
    0  & 0      &          &
    \\
    \tilde D   & 0  & \ddots   &
    \\
        & \ddots &  \ddots  & 0
    \\
        &        &   \tilde D  & 0
\end{bmatrix}
\]
Then using the above, for the AR-1 model
\[
\Big(G(\theta)x\Big)_\theta = -Dx
\]
And similarly
\[
\Big(G(\theta)^T\lambda\Big)_\theta = -D^\top\lambda
\]
We now have the expressions 
\[
f_\theta(\theta, \ol{y}) := -\ip{\ol{\lambda^p}, D\ol{x}}
\]
\[
f_{\theta \theta} = 0
\]
\[
\cH_{\theta}(\theta, \ol{y}) = -\begin{bmatrix} D^\top\ol{\lambda^p}\\ 0 \\ 0 \\ D\ol{x} \\ 0\end{bmatrix}
\]
When combined with the general results of section~\ref{sec:singular_ssm}, we get 

\begin{equation}
\label{eq:singular_ssm_der}
\begin{aligned}
v(\overline \theta) =& f(\overline \theta, \overline y)\\
v_\theta(\overline \theta) =& f_\theta(\overline \theta, \overline y) =  -\ip{\ol{\lambda^p}, D\ol{x}}\\
v_{\theta\theta}(\overline \theta) =&- \cH_{\theta}(\overline \theta, \overline y)^\top \cH_y(\overline \theta, \overline y)^{-1} \cH_{\theta}(\overline \theta, \overline y).
\end{aligned}
\end{equation}

\subsection{Fast Tracking of Unemployment Rates}
\begin{figure}[h!]
\begin{center}
\includegraphics[scale=.5]{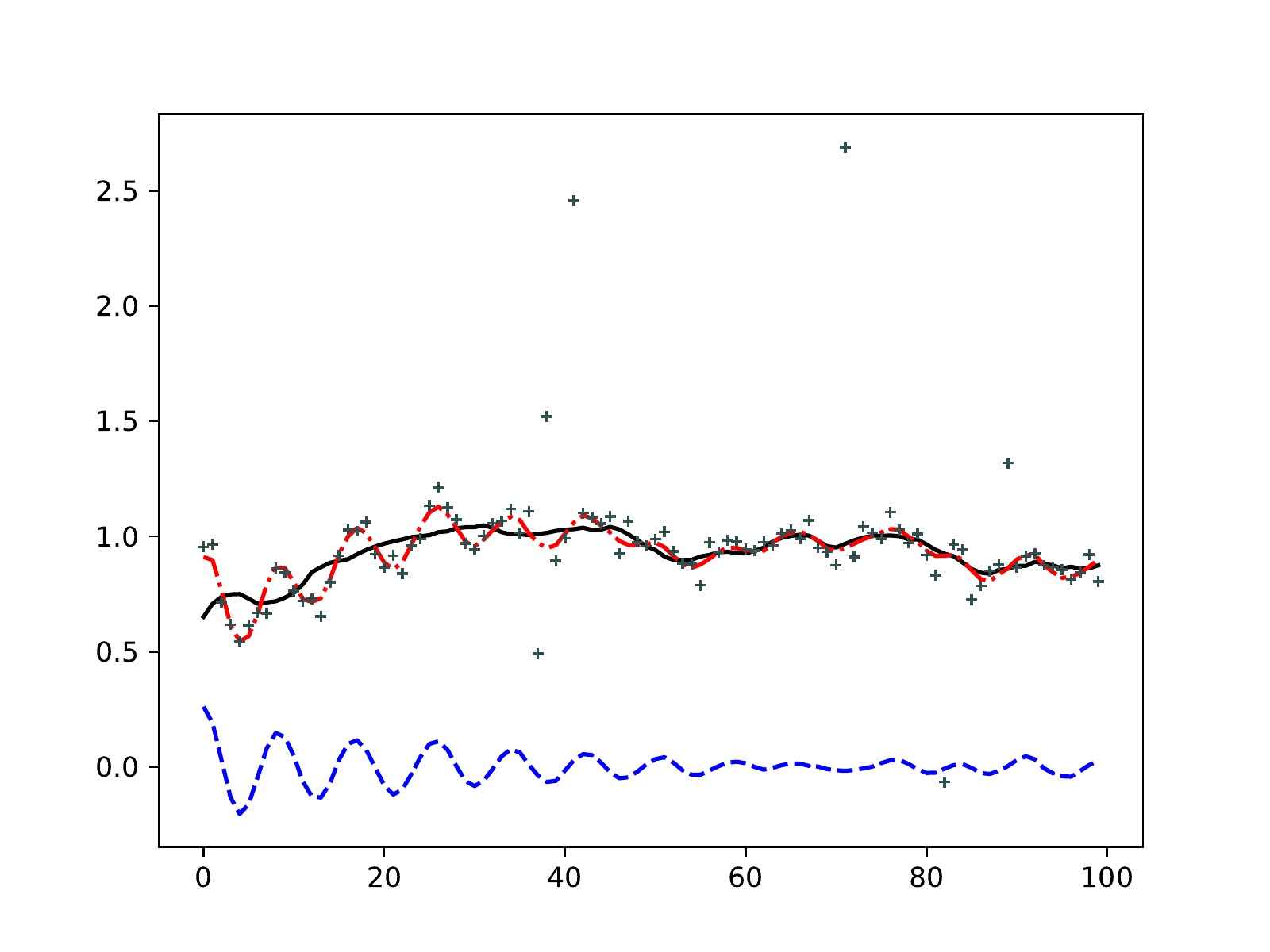}
\includegraphics[scale=.5]{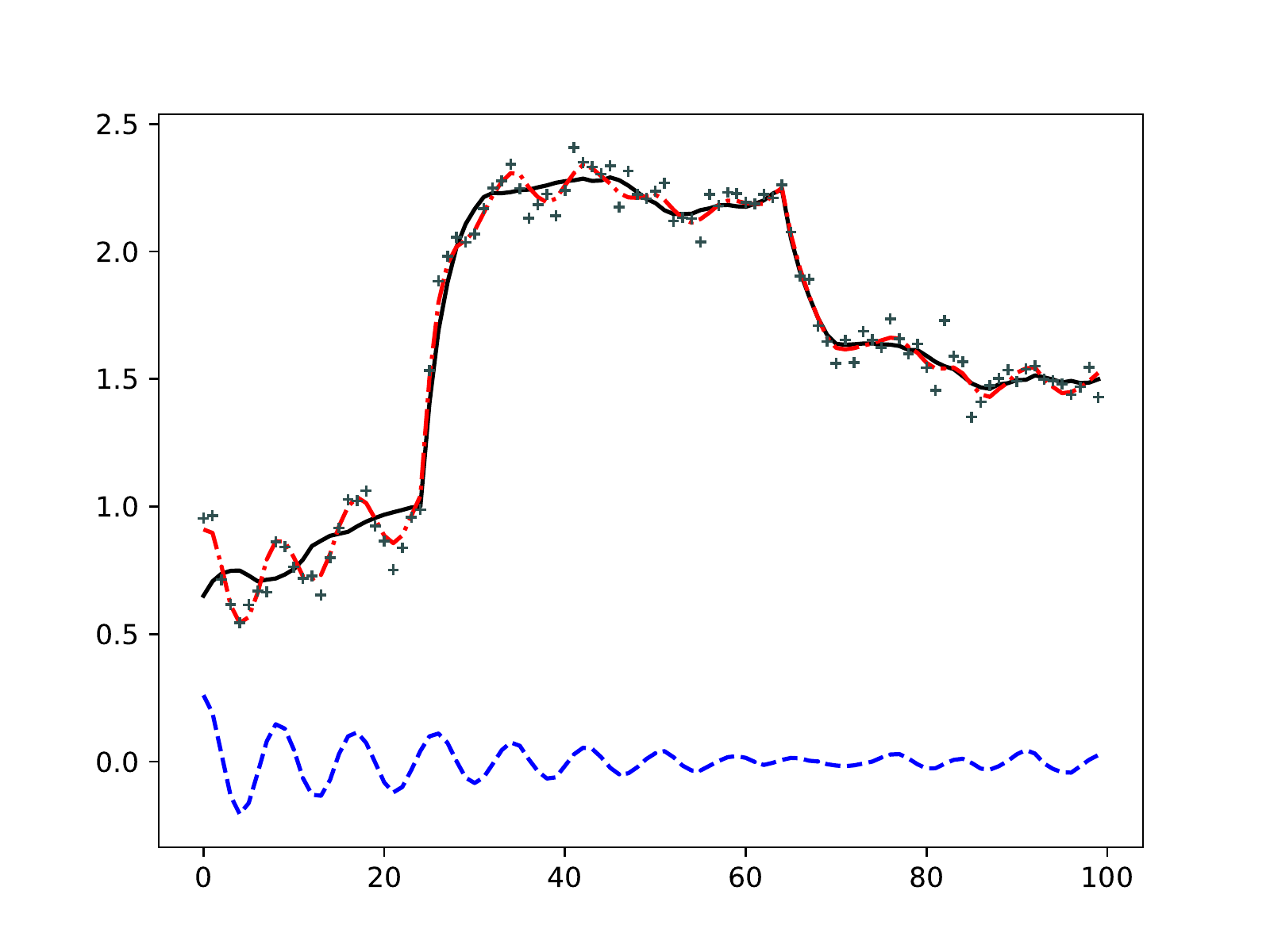}
\caption{\label{fig:data_ex}\small Example of generated data with outliers (top) and large jumps (bottom). Observations are shown in gray.}
\end{center}
\end{figure}

\begin{figure}[h!]
\begin{center}
\includegraphics[scale=.5]{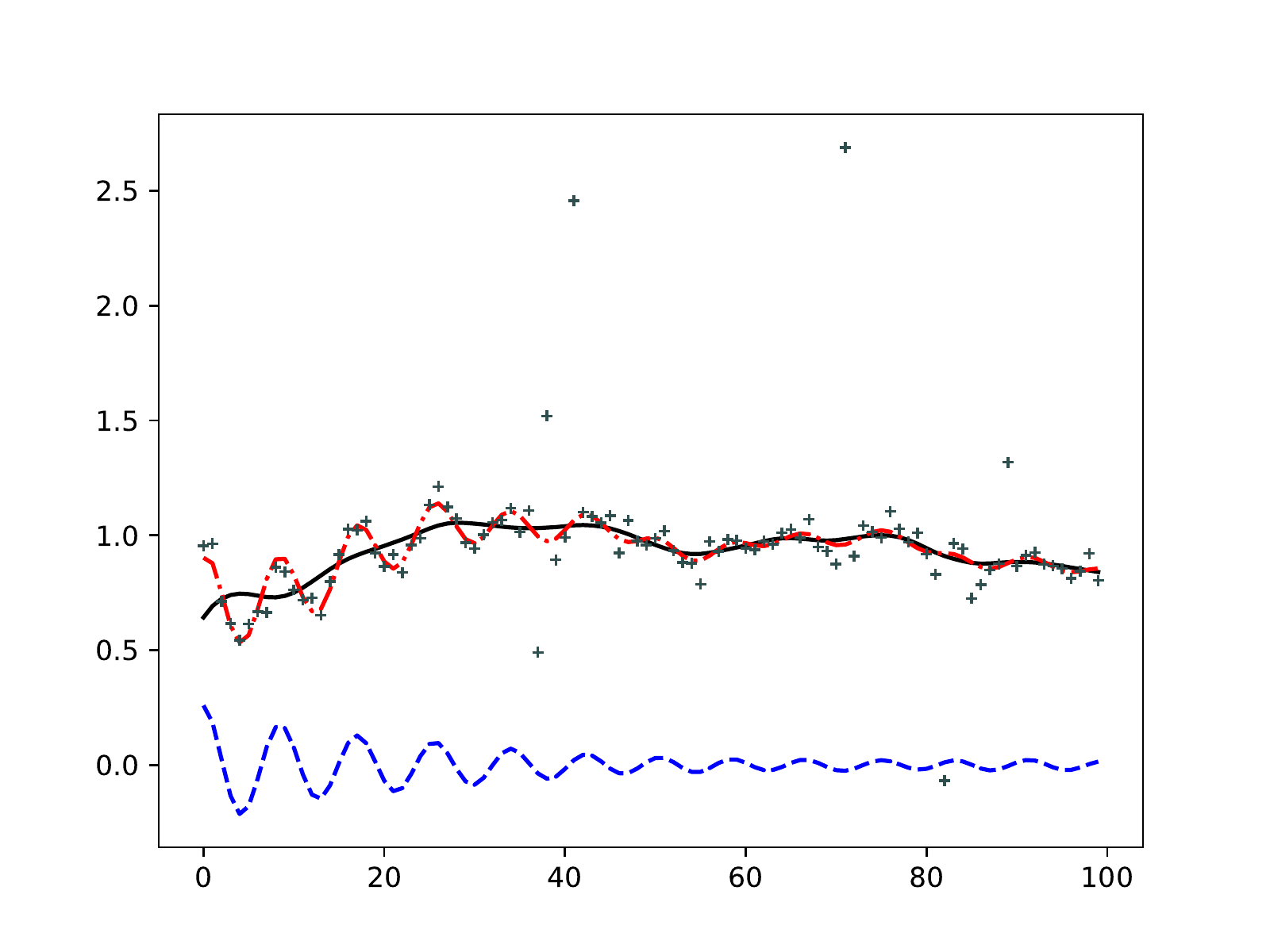}
\includegraphics[scale=.5]{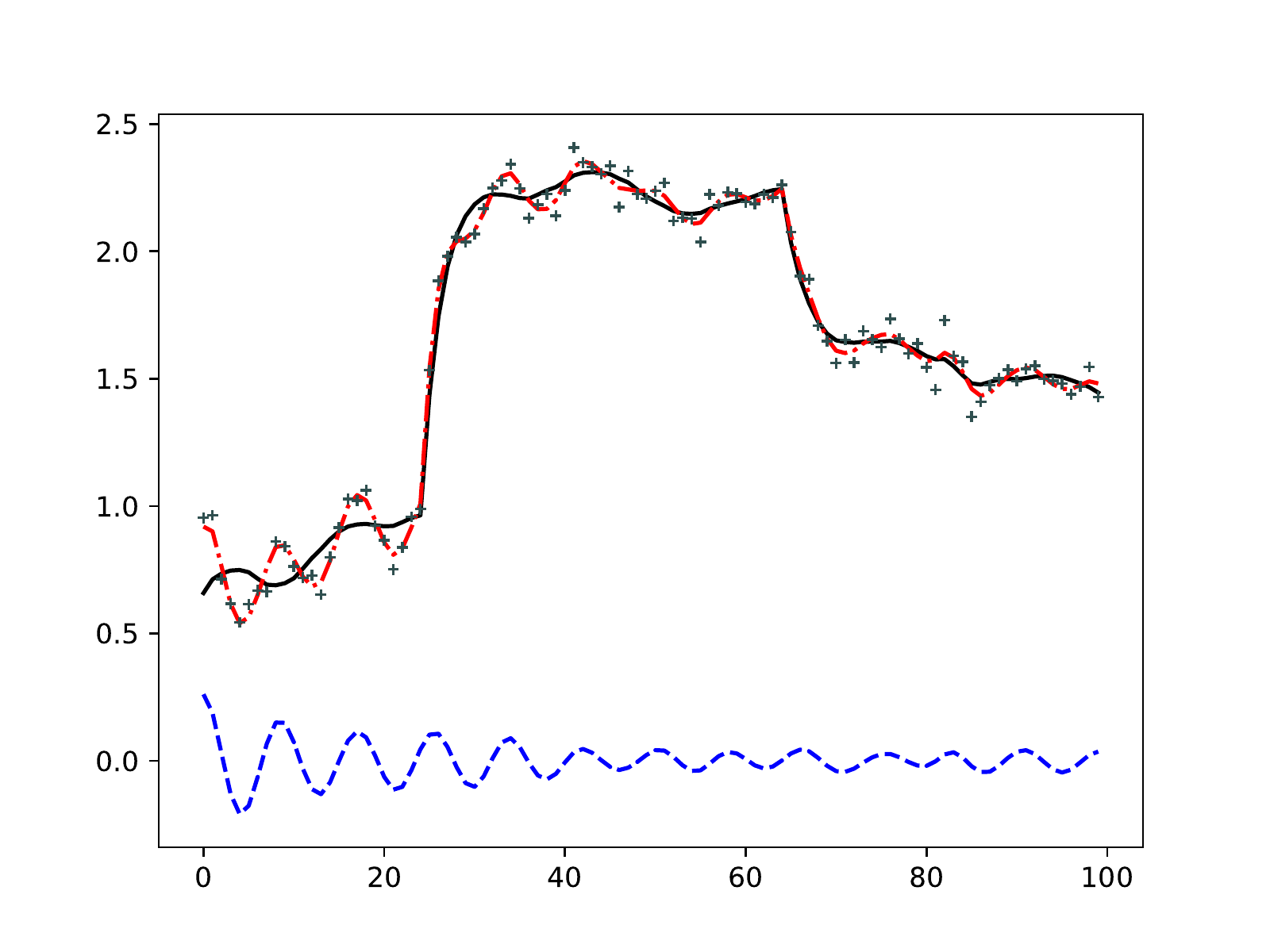}
\caption{\label{fig:data_est}\small Example of estimated solution run on data in Figure~\ref{fig:data_ex}. Top is ls/H run on data with outliers. Bottom panel shows T/ls run on data with large jumps added. 
Both are computed using Newton's method.}
\end{center}
\end{figure}

In this section, we apply the proposed approach to estimate parameters $(\ell_1, \ell_2, \gamma)$ 
for the structural unemployment model~\eqref{eq:state}---\eqref{eq:measurements}.
To test the approach, we generate ground truth parameters and then create synthetic data in order to compare model performance and speed using different formulations 
and algorithms.  
The data is generated by fixing parameters to reasonable values similar to those observed in practice, 
namely at $\begin{bmatrix}\ell_1 & \ell_2 & \gamma\end{bmatrix} = \begin{bmatrix}0.68 & 1.41 & -0.68\end{bmatrix}$,
 and applying the unemployment rate state space model 
\eqref{eq:state}---\eqref{eq:measurements} to generate the state as well as noisy observations. 
We then consider three cases: nominal errors, outliers in the observations, and jumps in the unemployment process. 
In the nominal cases we use variance parameters known to the smoother. 
To generate outliers we randomly select $10\%$ of measurements and add additional noise drawn from a $\mathcal{N}(0,1)$ Gaussian distribution. 
To generate large jumps we add large deviations  $.4, -.2$  at indices corresponding to $25$ and $65$. 
To Examples of the generated data are in Figure~\ref{fig:data_ex}. An example of the estimated using this data are in Figure~\ref{fig:data_est}.

\begin{table*}
\begin{center}
%\scriptsize
\small
\begin{tabular}{|c |c || c | c | c|| c | c | c|| c | c| c|}
\hline
& &\multicolumn{3}{c||}{Nominal} & \multicolumn{3}{c||}{Outliers} & \multicolumn{3}{c|}{Large Jumps}\\
\hline
&  & ls/ls & T/ls & H/ls & ls/ls & ls/T & ls/H & ls/ls & T/ls & H/ls \\
\hline
\multirow{4}*{Newton}&
	$||\hat{\theta}-\theta||^2$ & 0.158 & 0.143 & 0.1884& 1.72 & 0.115 & 0.142 & 0.296 & 0.032 & 0.099\\
	&Inner Iter & 19 & 250 & 206 & 14 & 71 & 102 & 21 & 329 & 373\\
	&Outer Iter & 9 & 5 & 10 & 6 & 15 & 13 & 13 & 9 & 14\\
	&Time & 11.5 & 36.5 & 29.7 & 6.3 & 20.2 & 19.3 & 13.3 & 48.1 & 49.5\\
\hline
\multirow{4}{*}{L-BFGS}&
	$||\hat{\theta} - \theta||^2$ & 0.158 & 0.137 & 0.184 & 1.72 & 0.115 & 0.141 & 0.296 & 0.034 & 0.099\\
	&Inner Iter & 19 & 228 & 207 & 19 & 86 & 88 & 32 & 375 & 519\\
	&Outer Iter & 10 & 18 & 10 & 9 & 17 & 12 & 16 & 14 & 15\\
	&Time & 6.2 & 33.9 & 26.5 & 6.3 & 19.4 & 14.0 & 10.6 & 51.4 & 67.5\\
\hline
\multirow{4}{*}{LM-Newton}&
	$||\hat{\theta} - \theta||^2$ & 0.158 & 0.069\footnote{Requires slightly better parameter initialization} & 0.184 & 1.04 & 0.12 & 0.142 & 0.296 & 0.028 & 0.099\\
	&Inner Iter & 16 & 385 & 180 & 32 & 70 & 114 & 24 & 302 & 357\\
	&Outer Iter & 12 & 15 & 20 & 25 & 28 & 19 & 18 & 15 & 21\\
	&Time & 7.2 & 51.6 & 26.4 & 14.4 & 20.0 & 21.0 & 11.3 & 43.2 & 48.1\\
\hline
\end{tabular}
\end{center}
\caption{\label{tab:tests} Table of results when run on generated data. The second row indicates the loss functions that were used where ls stands for least squares ($\frac{1}{2}||\dot||^2$), T stands for Students T with $\nu = 10$, and H stands for Hybrid with $\epsilon = .7$.}
\end{table*}


\label{sec:rates}

All algorithms are initialized at $\begin{bmatrix}0 & 0 & 0\end{bmatrix}$, 
except for LM-Newton on T/ls in the nominal case, which is initialized at $\begin{bmatrix}0 & 0 & 0.5\end{bmatrix}$,
as the standard zero initialization leads to bad results for this (nonconvex) case.  
 
In the first iteration the state is  initialized by propagating the initial $x_0$ through the dynamics for all time. 
In subsequent iterations the state is always initialized using the previous state solution. %previous solution is passed forward as initialized for the next state solve.
In all methods, the full state at each iteration is computed using Newton's method to find a saddle point of the augmented Lagrangian
\begin{equation}
\label{eq:auglag}
\begin{aligned}
AL(y,\theta) = & \ell^p (r^p)  + \ell^m(r^m) \\
&- \ip{\lambda^p, Q^{1/2}r^p - G(\theta)x + \zeta } \\
&-\ip{\lambda^m, R^{1/2}r^m - H(\theta)x + z } \\
& + \frac{1}{2}\|Q^{1/2}r^p-G(\theta)x+\zeta\|^2 \\
& + \frac{1}{2}\|R^{1/2}r^m-H(\theta)x+z\|^2
\end{aligned}
\end{equation}
using the Hessian of the Lagrangian provided in the Appendix. 

The value at an optimal point of (\ref{eq:auglag}) is the same as at an optimal point of (\ref{eq:Lag}). (\ref{eq:auglag}) is better conditioned which leads to faster convergence in practice.
For the outer iterations on the parameter space we compare a Newton method, L-BFGS, and a LM-Newton solver. The standard Newton and L-BFGS are from a standard python library. The LM-Newton solver is a quasi-Newton method where the Hessian is boosted by a parameter that is updated adaptively based on model performance. The results are summarized in Table~\ref{tab:tests}. % in the appendix. %An example of the estimates 

All methods work well for convex models. In the nonconvex case, the algorithms become more sensitive. 
In particular when $\ell^p$ is student's T, we have to boost $\ell^p_{rr}$ by a constant in order to make Newtons method converge. 
In practice this constant must be tuned depending on the parameter $\nu$ in the student's T function. 
Convergence is therefore sensitive to the choice of $\nu$ and boosting constant but a good rule of thumb is to choose $1 \leq \nu \leq 20$ 
and boost just enough to make the Hessian positive semidefinite.

\section{Discussion} 
We presented a general approach for parameter estimation in singular and non-singular Kalman smoothing models.
In particular we showed how to compute first and second derivatives 
of the value function (optimizing over state) with respect to the hidden parameters for both singular and nonsingular cases,
which captures a wide variety of models, including the motivating example. A simple numerical illustration
shows how the computed quantities can be used by a variety of optimization methods. The examples also show 
that when working with structural parameters, it pays off to have convex subproblems within each iteration 
of the value function. While non-convex losses such as Student's t are always appealing from a modeling perspective, 
when the problem is to find both the state and parameters, the resulting models are more fragile than those that use 
convex losses. This observation opens the way to future research in both theory and algorithm design.  

\smallskip

\noindent{\bf  Acknowledgment.}
We are very grateful to Jon Nielsen (Economic Council of the Labour Movement (ECLM)) for pointing us to the use of Kalman smoothers 
in structural unemployment models, and for teaching us about these models. 
Research of A. Aravkin and P. Zheng was supported by the Washington Research Foundation Data Science Professorship. 
J. Jonker was supported by the Boeing Data Science Research Grant. 
%Research of D. Drusvyatskiy was partially supported by the AFOSR YIP award FA9550-15-1-0237. The research of T. van Leeuwen was in part financially supported by the   Netherlands Organisation of Scientific Research (NWO) as part of research programme 613.009.032.

\bibliographystyle{abbrv}
\bibliography{kalmanSurvey,uwNav}

\newpage\newpage
\onecolumn
\appendix
\paragraph{Proof of Theorem~\ref{thm:invert}}
We reduce $\cH_{\cX}(\theta, \cX)$ in~\eqref{eq:HessLag} to block upper triangular form 
using invertible block row operations
\[
\begin{aligned}
\cR_1 & = (\ell_{rr}^p(r^p))^{-1} \cR_1 \\
\cR_2 & = (\ell_{rr}^m(r^m))^{-1} \cR_2\\
\cR_3 & = G(\theta)^{-T} \cR_3\\
\cR_4 & =  -G(\theta)^{-1}\left(\cR_4 - Q^{1/2}\cR_1 - Q^{1/2} (\ell_{rr}^p(r^p))^{-1}Q^{T/2} \cR_3\right) \\
\cR_5 & = \cR_5 - R^{1/2}\cR_2  - H\cR_4
\end{aligned}
\]
The resulting system is given by 
\[
\small
\begin{bmatrix}
I& 0 & - (\ell_{rr}^p)^{-1}Q^{T/2} & 0 &0 \\
0& I& 0 &0 & - (\ell_{rr}^m(r^m))^{-1}R^{T/2}\\
0 & 0 & I & 0& -G^{-1}H^\top \\
0 & 0 &0 & I & G^{-1}Q^{1/2} (\ell_{rr}^p(r^p))^{-1}Q^{T/2}G^{-T}H^\top \\
0 & 0 & 0 & 0 & \big(HG^{-1}Q^{1/2} (\ell_{rr}^p(r^p))^{-1}Q^{T/2}G^{-T}H^\top \\
 &&&& + R^{1/2} (\ell_{rr}^m(r^m))^{-1}R^{T/2}\big) 
\end{bmatrix}
\]
The invertibility of $\cH_\cX$ is thus equivalent to the invertibility of the symmetric positive semidefinite system~\eqref{eq:SchurLag}.

\paragraph{Hessian of the Lagrangian~\eqref{eq:auglag}}
\begin{equation}
AL_{yy}(y,\theta) = \begin{bmatrix}
G(\theta)^TG(\theta)+ H(\theta)^TH(\theta) & -G(\theta)^TQ^{1/2} & -H(\theta)^TR^{1/2} & G(\theta)^T & H(\theta)^T\\
-Q^{T/2}G(\theta) & \ell^p_{rr}(r^p) + Q^{T/2}Q^{1/2} & 0 & -Q^{T/2} & 0\\
-R^{T/2}H(\theta) & 0 & \ell^m_{rr}(r^m) + R^{T/2}R^{1/2} & 0 & -R^{T/2}\\
G(\theta) & -Q^{1/2} & 0 & 0 & 0\\
H(\theta) & 0 & -R^{1/2} & 0 & 0
\end{bmatrix}
\end{equation}

\end{document}